\documentclass[12pt]{article}

\usepackage{amsmath}
\usepackage{amssymb}
\usepackage{amsfonts}
\usepackage{amsthm}
\usepackage{mathrsfs}
\usepackage{url}
\usepackage[all]{xy}
\usepackage[pdftex]{graphicx}
\usepackage{fullpage}
\usepackage{indentfirst}
\usepackage{enumerate}
\usepackage{setspace}
\usepackage{epsfig}
\usepackage{fancyhdr}
\usepackage{hyperref}

  \newlength{\BiblioSpacing}
  \renewenvironment{thebibliography}[1]{%
    \begin{oldthebibliography}{#1}%
      \setlength{\parskip}{\BiblioSpacing}
      \setlength{\itemsep}{\BiblioSpacing}
  }%
  {%
    \end{oldthebibliography}%
  }

\newtheorem{thm}{Theorem}
\newtheorem{cnj}{Conjecture}
\newtheorem{lem}{Lemma}
\newtheorem{cor}{Corollary}[thm]         
 
\theoremstyle{definition}
\newtheorem{defn}{Definition}

\title{On the Classification of Universal Rotor-Routers}
\date{November 6, 2011}
\author{Xiaoyu He}

\begin{document}

\maketitle

\begin{abstract}
\noindent
The combinatorial theory of rotor-routers has connections with problems of statistical mechanics, graph theory, chaos theory, and computer science. A rotor-router network defines a deterministic walk on a digraph $G$ in which a particle walks from a source vertex until it reaches one of several target vertices. Motivated by recent results due to Giacaglia et al., we study rotor-router networks in which all non-target vertices have the same type. A rotor type $r$ is \textit{universal} if every hitting sequence can be achieved by a homogeneous rotor-router network consisting entirely of rotors of type $r$. We give a conjecture that completely classifies universal rotor types. Then, this problem is simplified by a theorem we call the Reduction Theorem that allows us to consider only two-state rotors. A rotor-router network called the compressor, because it tends to shorten rotor periods, is introduced along with an associated algorithm that determines the universality of almost all rotors. New rotor classes, including boppy rotors, balanced rotors, and BURD rotors, are defined to study this algorithm rigorously. Using the compressor the universality of new rotor classes is proved, and empirical computer results are presented to support our conclusions. Prior to these results, less than 100 of the roughly 260,000 possible two-state rotor types of length up to 17 were known to be universal, while the compressor algorithm proves the universality of all but 272 of these rotor types.
\end{abstract}

\section{Introduction}\label{intro}

The rotor-router model gives a deterministic analogue to random walks on directed graphs (digraphs). It is an offshoot of the abelian sandpile, or chip-firing, model introduced by Bak, Tang, and Wiesenfield \cite{BTW87} to study physical phenomena collectively known as ``self-organized criticality." Work by Dhar \cite{DHA06} and others on chip-firing gave new insights into phenomena such as the spreading of forest fires, the formation of fractal-like river networks, and other problems in statistical mechanics. Their work, notably the discovery of the sandpile group, laid the groundwork for further progress in chip-firing and rotor-routing. More recent combinatorial work on chip-firing can be found in the work of Diaconis and Fulton \cite{DF91}, and that of Levine \cite{LEV02}.

Rotor-routers, first introduced by Priezzhev et al. \cite{PDDK96} to study self-organized criticality, were then independently rediscovered several times in various connections. These include load-balancing, as in Rabani, Sinclair and Wanka \cite{RSW98}, and the Internal Diffusion Limited Aggregation (IDLA) model, as in Diaconis and Fulton \cite{DF91}. 

Rotor-routers have been studied extensively as a quasirandom analogue of random walks and Markov Chains \cite{HP09,LP07}, in particular for understanding IDLA. Recent progress shows that rotor-router aggregation satisfies similar, though much stronger, spherical asymptotics as random IDLA \cite{LP07}. Holroyd et al. \cite{HLMPPW} outlined a more general theory of rotor-routers on arbitrary digraphs, followed by the 2011 article \cite{GLPZ} by Giacaglia et al. that proved the periodicity of rotor-router hitting sequences and the conservation of palindromicity.

The universality problem, proposed by Propp, has been of interest since a solution can be viewed as a significant generalization of the aforementioned results by Giacaglia et al. Propp asks for the classification of rotor types, the fundamental units of rotor-router networks, that can model all others. These are called \textit{universal}. Any complete theory of rotor-routers requires an understanding of rotor universality, just as number theory requires an understanding of primes and ring theory an understanding of units.

Our main results towards this end are the Reduction Theorem of Section \ref{reductsect}, reducing the problem to the fairly manageable space of two-state rotors, and the introduction of the compressor algorithm in Section \ref{compressor}, which proves universality for a large number of rotors and significantly simplifies the problem for the rest. 

The relevant results from \cite{GLPZ} are presented in Section \ref{backmot} along with a rigorous description of the rotor-router model and the concept of universality. In Section \ref{examples}, we give two simple examples to elucidate our definitions, the latter of which is also essential to the problem at hand. In Section \ref{reductsect}, the Reduction Theorem is stated and proved. Section \ref{compressor} introduces the compressor, and its basic properties are proven. In Section \ref{empiricalresults}, we present various empirical results showing that the compressor algorithm almost always terminates in the rotor $12$ for unboppy rotors. Sections \ref{morebal} and \ref{abba} focus on the action of the compressor on specific classes of balanced rotors. Finally, Section \ref{future} makes explicit the many possible directions of future work hinted at in the rest of the report, and in Section \ref{acknowledgements} the people who made this work possible are acknowledged.

\section{Background and Motivation}\label{backmot}

On a finite directed graph $G = (V, E)$, allowing self-loops and multiple edges, a rotor-router network is imposed as follows. We choose a \textit{source} vertex as well as a set of \textit{target} vertices with outdegree $0$. Every non-target vertex must be able to reach at least one target. At each non-target vertex, we choose a \textit{rotor pattern} denoted $e_v ^{(1)}, e_v ^{(2)}, \ldots$ comprising an infinite periodic sequence of the directed edges emanating from $v$.

A particle travels through the network, beginning at the source and traveling until it reaches a target, whereupon it is returned to the source and commences a new walk. After its $n$-th visit to non-target vertex $v$, the particle leaves $v$ along edge $e_v ^{(n)}$. At this moment, we say that the first $n$ terms of the rotor pattern at $v$ have been \textit{used up}, or \textit{fired}. It is important to note that this scheduling policy is purely local and deterministic, and does not depend on the in-edge by which the particle enters $v$. The \textit{hitting sequence} of the network is defined as the infinite sequence of targets successively visited by this particle. Hitting sequences are eventually periodic; in \cite{GLPZ}, it is shown that they are also periodic from the start.

Each rotor pattern can be described by a partition $\pi$ of $\{1, 2, \ldots, n\}$ (this set is denoted $[n]$) where $n$ is its period, and $i, j \in [n]$ are partitioned into the same block of $\pi$ iff $e_v ^{(i)} = e_v ^{(j)}$. Each block of $\pi$ is associated with a distinct natural number, which gives a representation of the rotor that we call its \textit{rotor type}. There are an infinite number of such representations for every rotor pattern. For instance, if the pattern at a vertex $v$ with two out-edges $e$ and $e'$ is the sequence $e, e', e', e, e', e', \ldots$ with period 3, we can write its type as 1, 2, 2, or as 4, 7, 7. Typically if all the numbers are single digits, commas are omitted and we write 122 or 477 instead.

Often, we abuse these definitions and use the terms ``rotor pattern," ``rotor type," and ``rotor" interchangeably to avoid cumbersome language.

A rotor is said to have $m$ states if exactly $m$ distinct natural numbers appear in its type.  Two rotor types $r$, $r'$ are said to be \textit{equivalent} if they can represent the same rotor pattern. This is denoted $r \equiv r'$. Thus $122 \equiv 477$.

Given a rotor type $r$, the notation $|r|$ denotes its period length, and $r^{(k)}$ represents its $k$-th term, where $k$ is taken modulo $|r|$.

Propp introduced the idea that an entire rotor-router network can be thought of as a single meta-rotor with rotor pattern equivalent to its hitting sequence, since hitting sequences are always periodic.

To formalize this concept, we say that a rotor type $r$ can \textit{model} a rotor type $r'$ if there exists a \textit{homogeneous} rotor-router network, one containing only rotors of type $r$, with a hitting sequence equivalent to $r'$. Working off this idea, Propp defined a \textit{universal} rotor as a rotor that can model any other. Propp then raised the natural question: which rotor types $r$ are universal?

A rotor type is called \textit{palindromic} if it reads the same forwards and backwards (e.g. $121$), and \textit{block-repetitive} if it is the concatenation of uniform sequences of the same length $m$ (e.g. $1122$). The length $m$ is known as the \textit{block length} of such a rotor. We define:

\begin{defn}
By a \textit{block of length $m$} of a rotor type $r$, where $m$ divides $|r|$, we mean a subsequence of $r$ of the form $r^{(km+1)}$, $r^{(km+2)}$, ..., $r^{(km+m)}$, where $k \in \mathbb{N}$. We call this specific block the $(k+1)$-th block of length $m$ of $r$, denoted $r^{(k,m)}$.
\end{defn}

A block-repetitive rotor of block length $m$ is one such that every block of length $m$ is uniform. It is important to distinguish blocks from \textit{runs}, where the two are distinct in that $r^{(a+1)}, r^{(a+2)}, \ldots, r^{(a+b)}$ is always a run of $r$ but only a block if $b$ divides $a$. Uniform runs are important in Section \ref{morebal}.

It was proven by Giacaglia et al. \cite{GLPZ} that a network composed entirely of palindromic rotors has a palindromic hitting sequence, and similarly that a network composed entirely of block-repetitive rotors of block length $m\in \mathbb{N}$ has a block-repetitive hitting sequence with the same block length $m$. From these results it follows that neither palindromic nor block-repetitive rotor types can be universal. We believe these two classes encompass all nonuniversal rotors, so the following definition is natural:

\begin{defn}
A rotor type is \textit{boppy} if it is palindromic or block-repetitive. Otherwise, it is \textit{unboppy}.
\end{defn}

Henceforth, all bold uppercase letters or words will stand for classes of rotors. We let $\mathbf{ROT}$ denote the class of all rotor types, $\mathbf{PAL}$ the class of all palindromic rotors, $\mathbf{BR}$ the class of all block-repetitive rotors, $\mathbf{U}$ the class of all universal rotors, and $\mathbf{BOP} = \mathbf{BR} \cup \mathbf{PAL}$ the class of all boppy rotors. Also, $\overline{\mathbf{S}}$ denotes the complement of a set $\mathbf{S}$ of rotors (with respect to $\mathbf{ROT}$). In this notation, Giacaglia et al. \cite{GLPZ} showed that $\mathbf{U} \subseteq \overline{\mathbf{BOP}}$. We conjecture that more is true:

\begin{cnj}\label{cnjuniv}
As classes of rotor types $\mathbf{U} = \overline{\mathbf{BOP}}$.
\end{cnj}

Let the notation $\mathbf{C}(p, b)$, where $p \in \{0, 1\}$ and $b \in \mathbb{N}$, define the class of rotor types that are palindromic iff $p = 1$ and block-repetitive with block length $b$ (but no larger block length). The conjecture above can then be strengthened to:

\begin{cnj}\label{cnjstrong}
Every rotor type $r \in \mathbf{C}(p, b)$ is universal in $\mathbf{C}(p, b)$.
\end{cnj}

Conjecture \ref{cnjstrong} is only slightly stronger than Conjecture \ref{cnjuniv}, since Conjecture \ref{cnjuniv} implies that any rotor of class $\mathbf{C}(0,b)$ is universal in that class. It would only remain to show the same for the palindromic classes $\mathbf{C}(1, b)$. Finally, we also present a weakening that is more directly amenable to our methods:

\begin{cnj}\label{central}
Every rotor is universal in $\mathbf{C}(p, b)$ for some choice of $p$ and $b$.
\end{cnj}

In this paper, Conjectures \ref{cnjuniv} through \ref{central} will be our guiding focus, although we will make many related inquiries of independent interest.

The second half of this report focuses on the compressor configuration in Figure \ref{compfig}, so named because it tends to  shorten the lengths of the rotors upon which it is applied. The compressor can be viewed alternatively as a rotor-router network, as a mapping from rotor types to rotor types, or simply as a mapping from binary strings to binary strings. This last viewpoint is discussed extensively in Sections \ref{morebal} and \ref{abba}.

A related issue put into focus by the action of the compressor is the definition of a complexity measure for periodic binary strings. With the similar process known as the differential operator, or Ducci map, the complexity of such strings can be defined in terms of the length of the end cycle when the map is applied repeatedly \cite{KHO11}, or alternatively as the number of applications of the map before a term repeats. Since the compressor acts in an analogous way on balanced binary strings, we can also define a complexity measure of binary strings based on the compressor algorithm. This motivates our discussion of BURD and ab-ba rotors as the elements of fixed cycles of the compressor ``map" on binary strings. 

\section{Examples}\label{examples}

\begin{figure}[t]
  \centering
    \includegraphics[width=\textwidth]{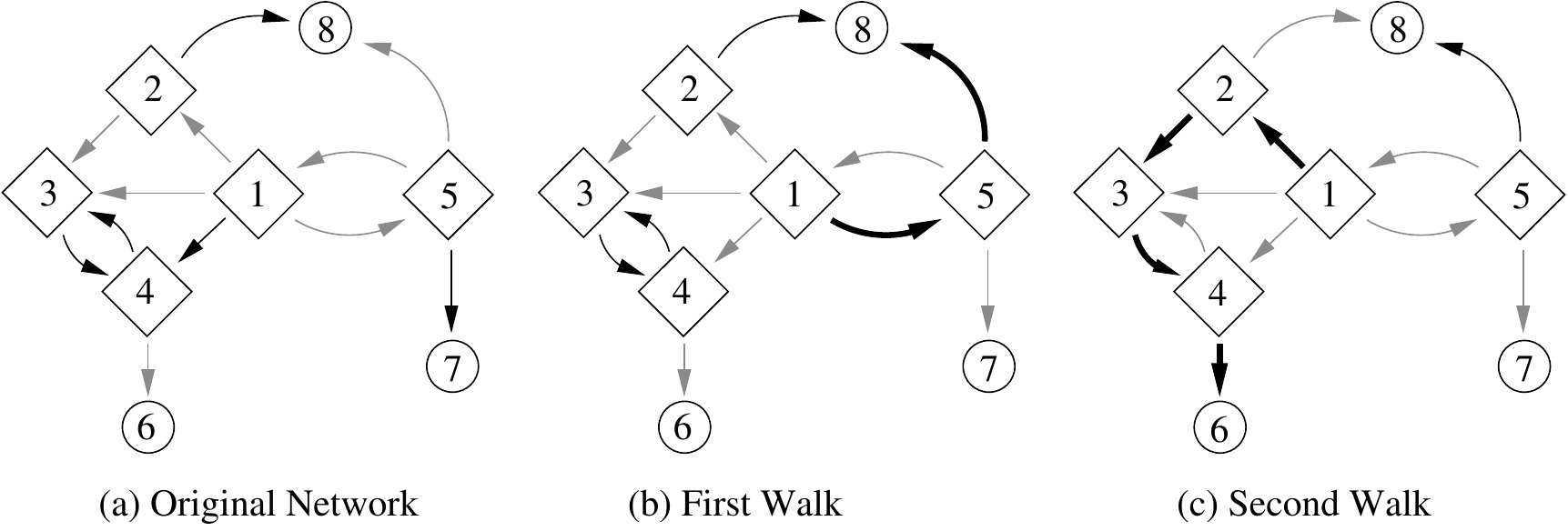}
  \caption{The first two walks in a simple rotor network. }
  \label{ex1}
\end{figure}

As a simple example of the rotor-router network, consider the directed graph in Figure \ref{ex1}. The underlying digraph, shown in Figure \ref{ex1}(a), consists of $5$ non-target vertices labeled $1, 2, 3, 4, 5$, of which vertex $1$ is the source, as well as three targets $6$, $7$ and $8$. The thin dark arrows denote the current state of the rotor at that vertex, or the direction by which the last particle left. In the other two diagrams, two walks of the particle to a target are traced by thicker dark arrows. Upon entering a new vertex, the particle shifts the rotor to the next possible exit direction (counterclockwise) and leaves along this new rotor setting. The rotor types are rather simple in this case, since the rotor patterns proceed counterclockwise; they are $1234, 12, 1, 12$, and $123$, at vertices $1, 2, 3, 4, 5$, respectively. The hitting sequence begins $8$, $6$, as shown.

From this example we already see some applications of rotor-routers. First, a rotor walk is the \textit{deterministic analogue of a random walk}, since rotor walks hit targets with the same frequencies as their associated random walks do in expectation \cite{HP09}. To be precise, rotor-routers are quasirandom in the sense that they are clearly nonrandom (unlike pseudorandom processes) and yet model random processes in some ways. Perhaps, if a randomized algorithm can be replaced with one involving rotor-routers, deterministic worst-case bounds may be given in place of probabilistic time bounds. It is not difficult to imagine how rotor-routers could be useful in load-balancing for parallel processing \cite{RSW98}. This ``derandomization" has also been used successfully to study the IDLA problem about random walks on the plane lattice \cite{LP07}.

Rotor walks share another property with random walks: the particle is \textit{not responsible for knowing its destination} (thanks to Yan Zhang of MIT for pointing this out). In conventional networking, each of the innumerable packets of information traveling through a network has some responsibility for finding its way to its destination. In a rotor-router network, however, particles are blindly shoved around but end up reaching targets in a predictable fashion. Thus rotor-routers could provide a robust mechanism for the distribution of information through large networks, especially if some regular pattern of delivery is required and the processing power at individual nodes is limited.

\begin{figure}[t]
  \centering
    \includegraphics[width=0.6\textwidth]{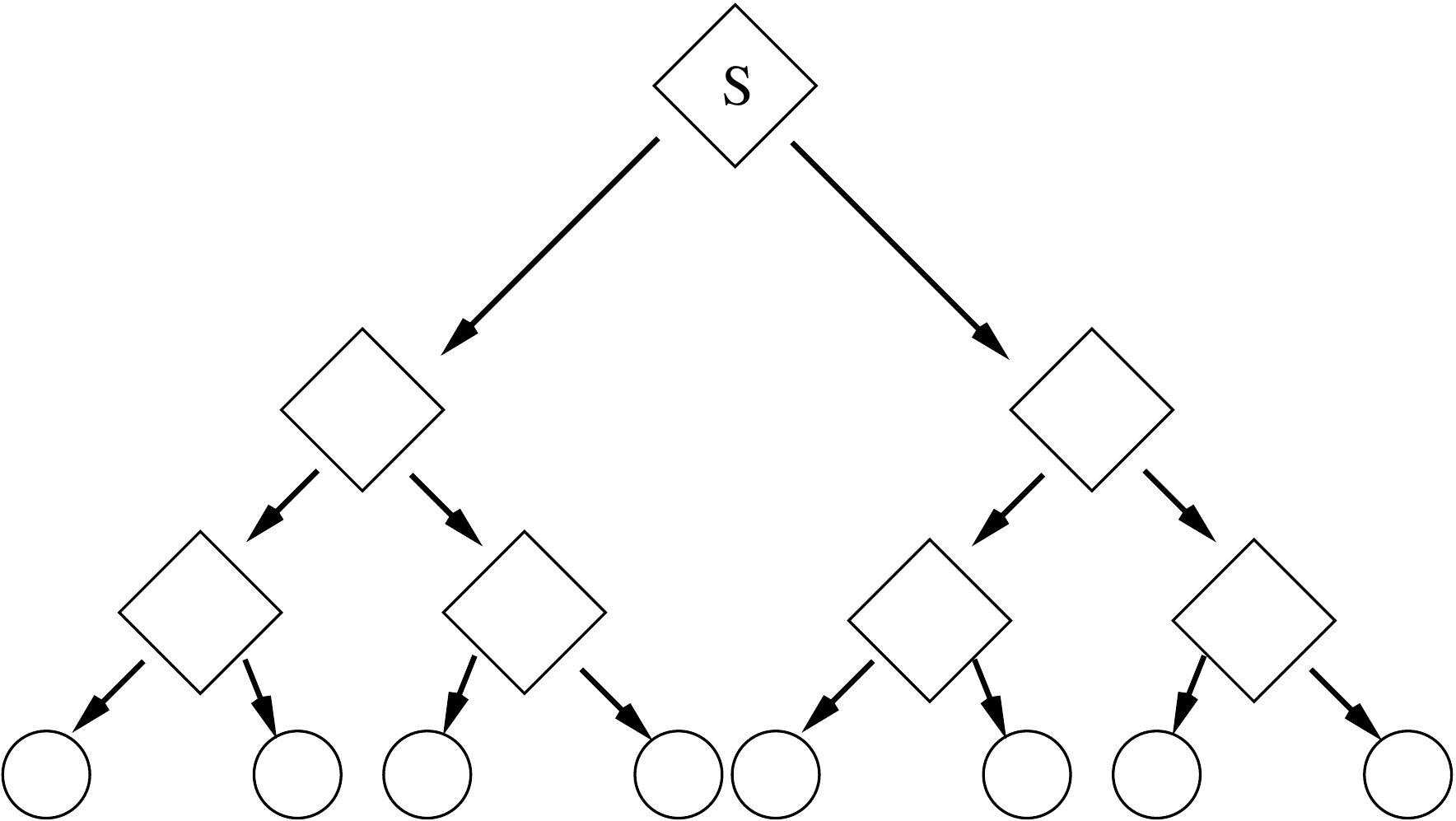}
  \caption{The binary tree network of depth $3$.}
  \label{ex2}
\end{figure}

We proceed to our second example, which is of direct interest to the problem at hand. We show that $12\in \mathbf{U}$. Consider the hitting sequence of a rotor network in the form of a complete binary tree of depth $n$ with $2^n - 1$ non-target nodes and $2^n$ targets, such that each rotor in the network is of type $12$. This means that particles leave each node in an alternating fashion. The source lies at the root of the tree. Figure \ref{ex2} shows such a tree of depth $3$. Now, we ask, what is the sequence of targets hit?

Let us show by induction that the hitting sequence is some permutation of the $2^n$ targets. In the case $n = 1$, we have a single non-target node, the source, so the two targets are hit in an alternating fashion. Now, suppose the statement is true for some $n \in \mathbb{N}$. Consider the binary tree of depth $n + 1$. It can be split into two parts: its first $n$ layers, forming a smaller binary tree, and its deepest layer, consisting of the $2^n$ non-target nodes directly preceding the targets. Using the inductive hypothesis, the first set can be replaced by a single source meta-rotor of type $1, 2, 3, \ldots, 2^n$, emitting particles to a permutation of the last layer of rotors. Now the result follows immediately by induction; this meta-rotor must emit two entire periods worth of particles to the last layer before any target is hit twice, since the rotors are of type $12$. It follows that all $2^{n+1}$ targets are hit once before any one is hit twice.

From the argument above, it follows that $12$ can make any rotor type of the form $1, 2, 3, \ldots, 2^n$ for $n \in \mathbb{N}$. In the spirit of Cauchy's proof of the AM-GM inequality, we modify the tree slightly to show that $12$ can make any rotor type of the form $1, 2, \ldots, m$, for any $m \in \mathbb{N}$. For each $m \geq 2$, choose $n$ such that $m \leq 2^n$. Consider the binary tree network of depth $n$, with $2^n$ targets, as before. Now, remove $2^n - m$ targets, routing their predecessors directly back to the root. This modification essentially ``hides" the extra targets the operation of the rest of the network is not affected (in the next section such a modification is called a \textit{destructive reduction}). The remainder of the hitting sequence has type $1, 2, \ldots, m$, as desired.

Finally, to form an arbitrary rotor type $a_1, a_2, \ldots, a_m$ of period $m$ from $12$, first construct the modified binary tree hitting $m$ targets, as above. Then, merge together the targets such that two targets $i, j\in [m]$ are combined iff $a_i = a_j$ (in the next section such an operation is called a \textit{merging reduction}). The hitting sequence is now the rotor type we want, so $12\in \mathbf{U}$.

Henceforth, all of our attacks on the universality problem will be based on the fact that $12$ is universal, since any rotor $r$ that can model $12$ is also universal.

\section{The Reduction Theorem}\label{reductsect}

The flexibility of the rotor-router model allows room for the generalization of results about two-state rotors to any number of states. In this section, we propose and prove a theorem regarding the reductions of multi-state rotors to two-state rotors.

A rotor type can often be thought of as a special case of another type with more states. The rotor type $12232$, for instance, is a special case, or \textit{reduction}, of the type $12434$, since we may direct each edge in a $12434$-type rotor numbered $4$ in the same direction as the one numbered $2$. This reduction is described informally as $4\rightarrow 2$, and called a \textit{merging reduction}. Along the same lines, $123$ is also a special case of the type $12434$, as each edge in a $12434$ rotor numbered $4$ can be turned into a self-loop, eliminating all the $4$'s. This reduction is denoted $x(4)$ informally, and called a \textit{destructive reduction}. In general we can always reduce a rotor by merging two of its states or deleting one altogether. The following lemma about reductions gets us close to our goal, Theorem \ref{reduction}.

\begin{lem}\label{reductlem}
If $r\in \overline{\mathbf{BOP}}$ has more than two states, and all of its merging reductions are boppy, then there exist three states of $r$, which we assume to be $1, 2, 3$, such that the reductions $3\rightarrow 1$ and $3 \rightarrow 2$ of $r$ are block-repetitive with coprime block lengths and $2\rightarrow 1$ is palindromic.

\begin{proof}
Because $r \in \overline{\mathbf{BOP}}$, there exists some integer $k$ such that $r^{(k)}$ and $r^{(|r|+1-k)}$ differ. Without loss of generality, call these two states $1$ and $2$, respectively. Consider another state, which we call $3$, without loss of generality. By our assumptions, if we apply $3\rightarrow 1$, the rotor becomes boppy. It cannot become palindromic, since $r^{(k)}$ and $r^{(|r|+1-k)}$ still differ, so it must become block-repetitive. Let its block size be $a$. Every block of $2$'s in $r$ must then have length divisible by $a$. Similarly, if we reduce $3\rightarrow 2$, the rotor becomes block-repetitive with some block size $b$, so the $1$'s must occur in blocks of length $b$. Now if $a$ and $b$ share any common factors $d$, both reductions are block-repetitive with block length $d$. Also, each block of $2$'s has length divisible by $d$, as does each block of $1$'s. Therefore, the remaining $3$'s must also be in blocks of length $d$, so that $r$ is itself block-repetitive with block length $d$. A contradiction has been reached so $a$ and $b$ must be coprime. 

Suppose these two merging reductions, $3\rightarrow 2$ and $3 \rightarrow 1$ are both block-repetitive, and consider the reduction $2\rightarrow 1$. First, suppose it too is block-repetitive, with block length $c$. It follows that the $3$'s occur in blocks of length $c$. Also, if $a$ or $b$ shares factors with $c$, then we have the same contradiction as before, that $r$ itself is block-repetitive. We have proven that $a$, $b$, and $c$ are pairwise coprime.

Choose any $k$ not a multiple of $abc$. Then, $r^{(k)}$ and $r^{(k+1)}$ lie in the same uniform block for at least two of the three numbers $a$, $b$, $c$. It follows that $r^{(k)} = r^{(k+1)}$. We have immediately that $r$ is block-repetitive with block length $abc$. Thus we have reached a contradiction and $2\rightarrow 1$ must give a palindromic reduction instead, as desired.
\end{proof}
\end{lem}

For the rotor type $111333332222111111222222222233333111$, five of its six possible reductions are boppy; two are block-repetitive and three are palindromic. Therefore, the following theorem may come as a surprise.

\begin{thm}[The Reduction Theorem]\label{reduction}
Every $r \in \overline{\mathbf{BOP}}$ with more than two states has an unboppy reduction.
\begin{proof}
Suppose there exists $r \in \overline{\mathbf{BOP}}$ such that all reductions of $r$ lie in $\mathbf{BOP}$.

Find some three states $1$, $2$, $3$ satisfying the conditions of Lemma \ref{reductlem}. Both reductions $3\rightarrow 1$ and $3\rightarrow 2$ of $r$ are block-repetitive with some coprime block lengths $a, b$, and so every state not $1$, $2$, or $3$ occurs in blocks of length $ab$. Therefore, we can discount them entirely; if we can show that the subsequence of $r$ containing only $1$, $2$, or $3$ has an unboppy reduction, then it is automatic that $r$ itself has an unboppy reduction by the same merging or destructive reduction rule. Indeed, if the reduction of this subsequence is not palindromic the corresponding reduction of $r$ cannot be palindromic. Similarly, if the reduction of this subsequence is not block-repetitive, then the corresponding reduction of $r$ just has some blocks of length $ab$ inserted, and is still not block-repetitive. Thus it suffices to assume that $r$ has exactly three states $1$, $2$ and $3$.

Let us prove that, in this case, the destructive reduction $x(3)$ is unboppy. The resulting type is not palindromic; otherwise, $r$ must have been palindromic, since $2\rightarrow 1$ is also palindromic. It remains to show that this rotor is not block-repetitive. If it is, the block length must be a multiple of $b$, since all $1$'s occur in such blocks, and a multiple of $a$, since all $2$'s occur in such blocks. Hence, the resulting type must have block length $ab$. But then, the $3$'s removed must also have been in blocks of $ab$ in order for the reductions $3\rightarrow 2$ and $3\rightarrow 1$ to have been block-repetitive. Hence $r$ itself was block-repetitive, which is a contradiction. Thus, $x(3)$ is unboppy. It follows that one of the reductions of $r$ lies in $\overline{\mathbf{BOP}}$.
\end{proof}
\end{thm}

After the Reduction Theorem, it suffices to consider rotor types of only two states; if Conjecture \ref{cnjuniv} is proven for every two-state rotor, then it follows in general by the Reduction Theorem. To tackle this reduced problem, we introduce the compressor in the next section.

\section{The Compressor}\label{compressor}
\begin{figure}[t]
  \centering
    \includegraphics[width=0.18\textwidth]{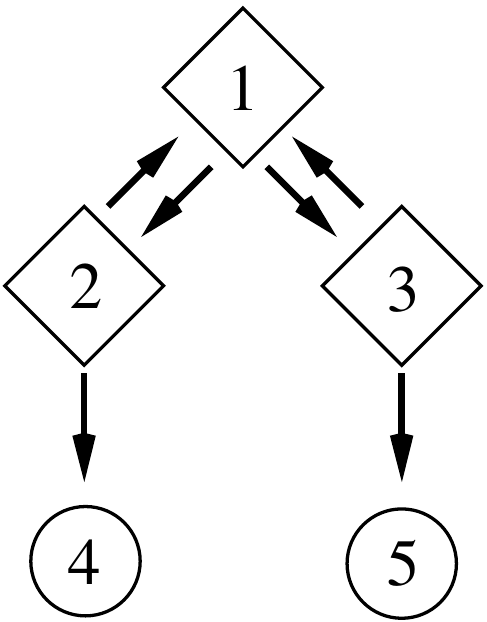}
  \caption{The compressor network.}
  \label{compfig}
\end{figure}

Here we introduce a rotor-router network of fundamental importance to the universality problem, the compressor. In one instance a special case of the compressor was used previously by Propp to prove the universality of $112$, but no recognition of its general significance was made. The rotor-router network in Figure \ref{compfig} is the compressor, a network that applies only to two-state rotors. It consists of three rotors of the same type $r$ as well as two targets. The source rotor is labeled $1$ and the targets are $4$ and $5$. By symmetry, there are four essentially different variations of this configuration, based on the order of the out-edges of rotors $2$ and $3$ (a type $12$ rotor at rotor $2$ could represent either $1, 4, 1, 4, \ldots$ or $4, 1, 4, 1, \ldots$).

Each configuration can be viewed as a mapping from the set of rotors to itself. The four variations of the compressor will be called $UU$, $UD$, $DU$, and $DD$. In each case, the first letter represents whether rotor $2$ points up to rotor $1$ first or down to rotor $4$, and the second letter represents whether rotor $3$ points up to rotor $1$ or down to $5$ first. For each rotor $r$, we define the mappings $UU, UD, DU, DD:\mathbf{ROT}\to \mathbf{ROT}$ to take $r$ to the hitting sequence of the corresponding compressor variation with three copies of $r$.

The compressor is a special case of the binary tree on three non-target vertices (with four targets), with two of the targets ``hidden," in the language of Section \ref{backmot}. In the language of the previous section, $UU(r)$, $UD(r)$, $DU(r)$, and $DD(r)$ are all destructive reductions of $BT(r)$, which will stand for the result of applying the binary tree to a rotor. Note that $BT$ takes two-state rotors to four-state rotors, and that $BT(r)$ can be formed from $UU(r)$ and $DD(r)$ by splitting these into smaller subsequences, and pasting these together; one of them will be given the labels $3, 4$ instead of $1, 2$. The same can be said of $UD$ and $DU$. This idea is helpful in Sections \ref{morebal} and \ref{abba}.

To illustrate this idea, consider the rotor $112212$, for instance. We find $UU(112212) = 454545 \equiv 121212$ and $DD(112221) = 445545 \equiv 334434$ (how this can be quickly determined will become clear presently), while $BT(112212) \equiv 334412132412$, an ``enmeshing" of $UU$ and $DD$.

To prepare for the next two results we define:

\begin{defn}
A \textit{balanced} rotor type is a two-state rotor type such that its two states occur with equal frequency. The class of all balanced rotors is denoted $\mathbf{BAL}$.
\end{defn}

The following result shows that $UD$ maps two-state rotors to balanced rotors. Henceforth, we assume without loss of generality that the two states of any two-state rotor are $1$ and $2$, and that the first term in the rotor is $1$.

\begin{thm} \label{UDbalanced}
If $r\in \mathbf{ROT}$ has two states, then $UD(r)\in \mathbf{BAL}$.

\begin{proof}
Suppose a period of $r$ has $a$ $1$'s and $b$ $2$'s. When exactly $|r| = a+b$ periods of the source rotor $1$ have been fired in network $UD$, exactly $a$ periods of rotor $2$ and exactly $b$ periods of rotor $3$ will have fired, so the graph is back to its original state. Thus some number of full periods of the hitting sequence has been reached. Since $ab$ of each target have been hit, $UD(r)\in \mathbf{BAL}$.
\end{proof}
\end{thm}

The proof of Theorem \ref{UDbalanced} also yields an $O(n^2)$ upper bound for the period length of $UD(r)$, where $n = |r|$, better than the trivial $O(n^3)$ bound.

\begin{cor}
If $r\in \mathbf{ROT}$ has two states, then $\displaystyle |UD(r)| \leq \frac{|r|^2}{2}$ and $|BT(r)| \leq |r|^2$.
\end{cor}

In general this argument extends to give $O(n^d)$ bounds for the hitting sequence lengths of homogeneous acyclic digraph (DAG) networks with size $m$ and longest chain of length $d$, a significant improvement over the exponential bounds for general rotor networks.

\begin{thm}
Let a rotor network be imposed on a digraph $G = (V,E)$ with $m$ non-target vertices and all rotors of type $r \in \mathbf{ROT}$. Suppose that every directed cycle in $G$ contains the source of the network, and $d$ is the largest number of non-target vertices in a simple directed path in $G$ beginning at source. Then the hitting sequence of $G$ has period length at most $|r|^d$.

\begin{proof}
We can assume that $G$ has no directed cycles. Otherwise, let the set of in-edges of the source $v_0 \in V$ be $E_0 \subset E$. Modify $G$ by adding one target vertex $t_0$. For each $(v, v_0) \in E_0$, we replace it with the edge $(v, t_0)$. There are no directed cycles in the new graph $G'$, and the period of the hitting sequence of $G$ is the destructive reduction $x(t_0)$, in the notation of the last section, of that of $G'$. Hence it is sufficient to prove the bound for acyclic networks.

It is only necessary to show that after the particle has been fired from the source $|r|^d$ times, it has passed through each non-target vertex a multiple of $|r|$ times. For each non-target vertex $v \in V$, we let $d(v)$ denote the length of the maximal path from $v_0$ to $v$. Then $d(v) \leq d - 1$. Inducting on $d(v)$, it is easy to show that the number of times the particle passes through $v$ is a multiple of $|r|^{d-d(v)}$. The result follows immediately.
\end{proof}
\end{thm} 

As an example of Theorem \ref{UDbalanced}, consider the action of $UD$ on the rotor type $112$. We write out (with foresight) some periods of the rotor patterns:

\begin{enumerate}[1:]
\setlength\itemindent{25pt}  
\item 223223223
\item 114114
\item 551
\end{enumerate}

Since we are applying $UD$, the first term in rotor $2$ is $1$ (``Up") while the first term in rotor $3$ is $5$ (``Down"). Since the number of $2$'s in three periods of the source rotor $1$ corresponds to the length of two periods of $2$, and the number of $3$'s to the length of one period of $3$, after everything displayed above is used up the rotor network has returned to its original position. Therefore the hitting sequence must contain exactly two $4$'s and two $5$'s.

Because of Theorem \ref{UDbalanced}, it suffices to study balanced rotors.

\begin{thm}\label{nonincreasing}
If $r \in \mathbf{BAL}$, then $|UU(r)| \leq |r|, |UD(r)| \leq |r|, |DU(r)| \leq |r|,$ and $|DD(r)| \leq |r|$.

\begin{proof}
If two periods of rotor $1$ have been fired, then exactly one period of rotor $2$ and one period of rotor $3$ will have fired, because $r$ is balanced. The rotors have all returned to their original states. Therefore, since exactly $|r|$ targets have been hit in this sequence, the period of the hitting sequence must divide $|r|$.
\end{proof}
\end{thm}

For instance, consider how $UU(121221)$ is formed. We write out two periods of the source and one period of rotors $2$ and $3$:

\begin{enumerate}[1:]
\setlength\itemindent{25pt} 
\item 232332232332
\item 141441
\item 151551
\end{enumerate}

Since $232332 \in \mathbf{BAL}$, we see that when every term listed above is used up, the network returns to its original state. Hence the hitting sequence is of length dividing $6$, just like the original rotor. In this particular instance the hitting sequence is $121212$, or $12$, of length $2$, and the compressor has proven that $121221 \in \mathbf{U}$. This example also illustrates why every variation of the compressor preserves balance, Corollary \ref{balclosedcor1}, and hints at Theorem \ref{one21} in Section \ref{abba}.

Repeated compressor applications never increase a balance rotor's period, since the compressor preserves balance. Many rotors can be ``compressed" to rotors of short period, or even to $12$, as above. Explicitly, the compressor algorithm runs as follows. Starting with an arbitrary two-state rotor $r$, apply $UD$ to create balance. Thereafter, apply an arbitrary compressor variation to the result (in which we replace $4$'s by $1$'s and $5$'s by $2$'s) and repeat. In the vast majority of cases if $r$ is unboppy this algorithm ends in $12$, proving the universality of $r$. This is the first, and crudest, instance of the use of monovariants associated with the compressor. As will be seen in Sections \ref{morebal} and \ref{abba}, more monovariants further restrict the action of the compressor.

Finally, we sketch an algebraic way to study the compressor. For a given $r\in \mathbf{BAL}$, write it as an infinite periodic sequence $a_1, a_2, \ldots$ of period $2n$. Define $f(m)$ to be the position of the $m$-th $1$ in this sequence, and $g(m)$ to be the position of the $m$-th $2$. Also, define $F(m)$ to be the number of $1$'s in the sequence up through $a_m$, and $G(m)$ the number of $2$'s up through $a_m$. Let us find the position of the $m$-th $4$ in $UD(r)$ in terms of $f$, $g$, $F$, and $G$.

\begin{thm}\label{positionUD}
The $m$-th $4$ in $UD(r)$ occurs at position $F(G(f(g(m)))) + m$.
\begin{proof}
Consider the position of the $m$-th $4$ in $UD(r)$. It must come from directly the $m$-th $4$ in rotor $2$ of the compressor configuration. We count the number of $5$'s hit before this $4$ to find its absolute position in $UD(r)$. It came from the $m$-th $4$ in rotor $2$, at position $g(m)$ in that rotor. This was hit by the particle coming from the $g(m)$-th $2$ in the source rotor, which is at position $f(g(m))$, since $2$'s in rotor $1$ correspond to $1$'s in $r$. Now, exactly $G(f(g(m)))$ of the terms in the source up to this point are $3$'s, and they all went to rotor $3$. Of these, exactly $F(G(f(g(m))))$ hit a $5$ in that rotor. It follows that the $m$-th $4$ in $UD(r)$ occurs at position $m + F(G(f(g(m))))$.
\end{proof}
\end{thm}

Theorem \ref{positionUD} and its analogues for the other compressor variations give a number of results supporting our intuition. We give only one example here. Suppose $UD(r)\in \mathbf{PAL}$ for some $r$. Note that since the $m$-th and $(n - m + 1)$-th $4$'s in $UD(r)$ must be symmetric (i.e. their positions must add to $2n + 1$), it follows that

\begin{equation}
F(G(f(g(m)))) + F(G(f(g(n-m+1)))) = n.
\label{PAL1a}
\end{equation}

By the exact same manipulations on the $5$'s in $UD(r)$, we can show that

\begin{equation}
G(F(g(f(m)))) + G(F(g(f(n-m+1)))) = n.
\label{PAL1b}
\end{equation}

These equations say that the original rotor $r$ was ``almost palindromic." For instance:

\begin{thm}
If $r\in \mathbf{BAL}$, then $UD(r)\in \mathbf{PAL}$ only if $r$ starts and ends with the same term.

\begin{proof}
With equations (\ref{PAL1a}) and (\ref{PAL1b}) in mind, suppose for the sake of contradiction that $r$ starts with $1$ and ends with $2$ (without loss of generality). Find the smallest $k$ such that $a_k = 2$, so $a_i = 1$ for all $i < k$. Then, letting $m = 1$ in (\ref{PAL1a}), we see that

\begin{equation}
F(G(f(k))) + F(G(f(2n))) = n.
\label{Cons1}
\end{equation}

Next, note that $G(f(i)) = f(i) - i$, since the number of $2$'s up to the $i$-th $1$ is exactly the position of the $i$-th $1$ minus $i$, the number of $1$'s up to that point. Thus equation (\ref{Cons1}) simplifies further to

\begin{equation}
F(f(k) - k) + F(f(2n) - 2n) = n.
\label{Cons2}
\end{equation}

Furthermore, since $f(2n)$ is the position of the $2n$-th $1$, it is the last $1$ in the second period of $r$. Because a period of $r$ has length exactly $2n$, it follows that $f(2n) - 2n = f(n)$, the position of the $n$-th $1$. Also $F(f(n)) = n$, so now it remains to show that

\begin{equation}
F(f(k) - k) = 0
\label{Cons3}
\end{equation}

\noindent is impossible. Note that $f(k) \geq k$ trivially, and $f(k) \neq k$ because $a_k = 2$. Therefore, $f(k) - k \geq 1$, from which it follows that $F(f(k) - k) \geq 1$, since $a_1 = 1$. Thus (\ref{Cons3}) is impossible and a contradiction has been reached. The rotor type $r$ must start and end with the same term.
\end{proof}
\end{thm}

Similar local results can be deduced from Theorem \ref{positionUD} for block-repetitive rotors. Both empirically and heuristically from the type of argument above, the compressor only immediately fails to prove universality for unboppy rotors when they are close to boppy rotors (a few terms away). We turn to the empirical data in the next section.

\section{Computer Data and Universal Rotor Types}\label{empiricalresults}

Of all nontrivial two-state rotors, it was previously known that the rotor types $1^k 2$ (which denotes $k$ $1$'s followed by a $2$) were universal. By judicious application of the compressor, the following general classes have been proven universal: those of the form $(11)^k 21$, and those of the form $(12)^k 21 (12)^l$. A short proof of the universality of the second class is given in Section \ref{abba}.

In addition to these classes, we have proven through C programs the universality of the vast majority of all unboppy two-state rotor types with period up through 17. The rotor-router model, and in particular our compressor algorithm, lends itself to computer simulation. Here is a brief sketch of the programs written to generate our data.

Based on an input parameter $n$, the length of the rotors to consider, all two-state rotors of that length are generated. Those that are trivial, such as $11111$, or necessarily nonuniversal (i.e. boppy) such as $1221$, are immediately discarded. For each of the remaining rotors, the program applies $UD$ to the rotor to make it balanced, after which it repeatedly applies a pseudorandomly chosen compressor variation for up to $200$ times. If $12$ is reached within this subroutine, the loop terminates and the rotor is universal. Otherwise, up to $50$ more such random sequences of applications of the compressor occur, restarting each time with the original rotor. If at the end of this loop the rotor still has not been proven universal, it is recorded in an output file as undecided. The program outputs at termination the total number of undecided rotors found. The running time is necessarily (at least) exponential and proved prohibitive for values of $n$ greater than $17$.

\begin{table}[ht]
\singlespacing
  \begin{center}
  \begin{tabular}{|l|l|l|}
  \hline
  Length & Undecided & Fraction of Total \\
  \hline
  2 & 0 & 0.0000 \\
  3 & 0 & 0.0000 \\
  4 & 0 & 0.0000 \\
  5 & 0 & 0.0000 \\
  6 & 2 & 0.0625 \\
  7 & 0 & 0.0000 \\
  8 & 5 & 0.0390 \\
  9 & 9 & 0.0351 \\
 10 & 22 & 0.0430 \\
 11 & 0 & 0.0000 \\
 12 & 36 & 0.0088 \\
 13 & 0 & 0.0000 \\
 14 & 12 & 0.0015 \\
 15 & 52 & 0.0032 \\
 16 & 136 & 0.0042 \\
 17 & 0 & 0.0000 \\
  \hline
  \end{tabular}
  \end{center}
\caption{The number of unboppy rotors not proven universal.}
 \label{resulttab}
\doublespacing
\end{table}

In Table \ref{resulttab}, we summarize these results. For each rotor length up through $17$, the table gives the number of unboppy rotor types not proven universal by the compressor, and its fraction out of the total number of unboppy rotors of that length. In total, at most 272 are not universal out of roughly 260000. It seems plausible from this data that \textit{almost all} unboppy rotor types can be proven universal using the compressor alone, in the sense that the fraction of undecided rotors approaches $0$.

It should be stressed that these undecided rotors are not necessarily counterexamples to Conjecture \ref{cnjuniv}; for instance, the $7$ potential counterexamples of length $6$ or $8$ have been proven universal by first applying an auxilliary rotor-router network and then using the compressor algorithm. Although this method times out for longer rotors it shows that the occasional failure of the compressor algorithm is probably not due to the existence of a third class of nonuniversals.

From this data we conjecture that all unboppy rotors of prime length can be proven universal by the compressor. It is true for all the primes up through 17. 

Variations on the program above produced further data supporting our conjectures. Every balanced, unboppy rotor with period less than 24 is universal, and at most 9 such rotors with period 24 are not. In addition, every ab-ba rotor of length $2 \pmod{4}$ with period up through 42 is universal. Finally, we note that every rotor tested, up to length 17, has been shown universal in some class $(p, b)$, strongly supporting Conjecture \ref{central}.

\section{More on Balanced Rotors}\label{morebal}

We extend the theory of balanced rotors and prove some finer claims than those of Section \ref{compressor}, culminating in the definition of BURD rotors. In this section it is possible to view rotors almost entirely as binary strings upon which the compressor family of maps are applied. We start by refining the original notion of rotor balance.

\begin{defn}
An \textit{$n$-balanced} rotor is a two-state rotor type with fundamental period divisible by $n$ such that for each $k \in \mathbb{N}$, the block $r^{(k,n)}$ has equal numbers of $1$'s and $2$'s. The class of $n$-balanced rotors is denoted $\mathbf{n}$-$\mathbf{BAL}$. If a rotor $r$ is $2$-balanced, we also call it an \textit{ab-ba rotor}; the class $\mathbf{2}$-$\mathbf{BAL}$ of $2$-balanced rotors is also called $\mathbf{ABBA}$.
\end{defn}

Note that the compressor preserves $n$-balance. As an illustration, consider $DD(12211122)$. In this case the compressor looks like:

\begin{enumerate}[1:]
\setlength\itemindent{25pt} 
\item 2332223323322233
\item 41144411
\item 51155511
\end{enumerate}

Because $12211122\in \mathbf{4}$-$\mathbf{BAL}$, we can see that when $23322233$ is fired in rotor $1$, balanced blocks $4114$ and $5115$ are used up in rotors $2$ and $3$. The next block of length eight, which is again $23322233$, corresponds to balanced blocks $4411$ and $5511$ in rotors $2$ and $3$. After each balanced block of eight is fired from the source, the corresponding two blocks of four are fired in rotors $2$ and $3$, producing the desired $4$-balanced effect. Formally, we have:

\begin{thm}\label{balclosed}
For any $n\in \mathbb{N}$, the class $\mathbf{2n}$-$\mathbf{BAL}$ is closed under $UU, UD, DU,$ and $DD$.
\begin{proof}
For $UD$ and $DU$, Theorem \ref{nequiv} gives a stronger result. We need only prove the result for $UU$ and $DD$. Consider $UU(r)$ for an $2n$-balanced rotor $r$. We wish to show that for each $k \in \mathbb{N}$, the $k$-th block of length $2n$ in $UU(r)$ is balanced. When the block $r^{(k, 4n)}$ is fired from the source, half of the particles enter rotor $2$, and the other half enter $3$. From $2$ and $3$, each fires according to the block $r^{(k,2n)}$, so we hit the two targets $n$ times each. Hence this block of $2n$ of the hitting sequence is balanced, as is every other by the same argument.
\end{proof}
\end{thm}

The two most important special cases are:

\begin{cor}\label{balclosedcor1}
The classes $\mathbf{BAL}$ and $\mathbf{ABBA}$ are closed under $UU, UD, DU,$ and $DD$.
\end{cor}

For the following results it is natural to define a coarser type of equivalence, where only the frequencies of labels within blocks are taken into consideration.

\begin{defn}
Two $r, r'\in \mathbf{ROT}$ are \textit{n-equivalent} if there exists $s\in \mathbf{ROT}$ equivalent to $r$ with the same states as $r'$ such that for each $k \in \mathbb{N}$, as multisets $r'^{(k,n)} = s^{(k,n)}$. We write $r \equiv_n r'$. In particular, if $r \equiv_1 r'$, then $r \equiv r'$.
\end{defn}

Now that we have the requisite definitions, we can state and prove Theorem \ref{nequiv}.

\begin{thm} \label{nequiv}
If $r\in \mathbf{2n}$-$\mathbf{BAL}$, then $UD(r) \equiv_n DU(r) \equiv_n r$.

\begin{proof}
We show $UD(r) \equiv_n r$. Suppose by induction that this has been shown for the first $k$ blocks of length $n$ in $UD(r)$ and $r$. Then, consider $r^{(k,2n)}$. At the point that the source rotor first enters this block of its cycle, by the balanced property of the previous blocks the two other rotors have just used up their first $k$ blocks of length $n$, respectively. Now, in the next $2n$ firings of the source, the next block of $n$ of rotor $2$ is used up. There are the same number of $4$'s in this block as there were $2$'s in $r^{(k,n)}$. The corresponding block of rotor $3$ also has the same number of $5$'s as there are $1$'s in $r^{(k,n)}$, so a block of length $n$ is completed and we are done by induction.
\end{proof}
\end{thm}

In particular, when Theorem \ref{nequiv} is applied to $2$-balanced, or ab-ba rotors, it shows that all ab-ba rotors are fixed points of $UD$ and $DU$.

\begin{cor}\label{abbafixed}
If $r\in \mathbf{ABBA}$, then $UD(r) \equiv DU(r) \equiv r$.
\end{cor}

Theorem \ref{nequiv} exhibits a possible obstacle of the compressor algorithm. If a rotor $r\in \mathbf{4n}$-$\mathbf{BAL}$ but $r \not\in \mathbf{2n}$-$\mathbf{BAL}$ for some $n$, then $UD$ and $DU$ will never make it $2n$-balanced; some blocks of length $4n$ become permanently skewed in one direction. For instance, $UD(11121222)=54554454$, and the imbalanced number of $1$'s and $2$'s in the first block of length $4$ is not resolved.

However, an almost exactly opposite result holds for $UU$ and $DD$, showing that they distribute the states more evenly. Recall that a balanced run of a rotor $r$ is a subsequence $r^{(a + 1)}, r^{(a+2)}, \ldots, r^{(a+b)}$ such that the same number of $1$'s appear as $2$'s. For the purposes of the next theorem, it is simpler to treat the rotor-router particle as multiple reincarnations, where it becomes a new, distinct particle whenever it returns to the source. Thus the particle that uses up the $k$-th term of rotor $1$ in the compressor is the $k$-th particle. In $UU(r)$ and $DD(r)$ there are two runs that together correspond to the original run in $r$ in the sense that the targets in those subsequences were hit by exactly those particles using up that run in the source.

\begin{thm} \label{UUDDdistribute2}
If $r\in \mathbf{BAL}$ is such that the runs $r^{(1)}, r^{(2)}, \ldots, r^{(2a)}$ and $r^{(2a + 1)}, \ldots, r^{(2a+2b)}$ are balanced, then the two runs in $UU(r)$ and $DD(r)$ corresponding to the $(2a+1)$-th through $(2a+2b)$-th particles to leave the source are both balanced and have lengths that add to $2b$. This is true unless the run $r^{(a + 1)}, \ldots, r^{(a + b)}$ is constant, in which case exactly one of the balanced runs exists and has length exactly $2b$.
\begin{proof}
Consider first $UU(r)$. The condition that $r{(1)}, r^{(2)}, \ldots, r^{(2a)}$ is balanced simply tells us that the same number of terms, namely $a$, of rotors $2$ and $3$ have been fired. The second condition shows that exactly $b$ particles among the next $2b$ enter rotor $2$, and the rest enter rotor $3$. Consider the number of $4$'s in rotor $2$, which is equal to the number of $5$'s in rotor $3$, coming from the run $r^{(a + 1)}, \ldots, r^{(a + b)}$; whatever this number is, it is the number of $4$'s and $5$'s in the corresponding balanced run of $UU$, and $2b$ minus this number is the corresponding size of the balanced run in $DD$, by the same logic. One of these numbers is zero iff $r^{(a + 1)}, \ldots, r^{(a + b)}$ is constant, so we are done.
\end{proof}
\end{thm}

Unless the exception occurs, we say that the run in $r$ is \textit{split up} between $UU(r)$ and $DD(r)$. Note that two periods of $r$ correspond to one period of $UU(r)$ and one of $DD(r)$, so an individual run can be split up twice, independently.

For an example of Theorem \ref{UUDDdistribute2}, consider again the rotor $12211122$, which has $UU(12211122) = 54454455$ and $DD(12211122) = 45454554$. We say that the first occurrence of the balanced run $1122$ is split up between the $45$ in positions three and four of $UU$ and the corresponding $45$ of $DD$. The second occurrence of $1122$, meanwhile, is not split up because the run at half the position with half the length is the last $22$ in the rotor, and is uniform. Instead, this balanced run corresponds to the $4455$ that ends $UU$.

These ideas about balanced and uniform runs motivate two more notions.

\begin{defn}
Define a \textit{balanced run decomposition} (BRD) of a balanced rotor $r$ as a partition of the period of $r$ into contiguous balanced runs. Define the balanced run type of this BRD as the infinite periodic sequence $(b_1, b_2, \ldots)$ such that $2b_i$ is the length of the $i$-th balanced run, modulo the size of the partition. Let the \textit{balance coefficient} of $r$, denoted $b(r)$, be the maximal number of runs in any BRD of $r$, divided by $|r|$.
\end{defn}

A rotor may have multiple BRD's if we choose to merge adjacent balanced runs. Thus $122112$ can have BRD $12|21|12$ or $12|2112$ or $122112$. Note that Theorem \ref{UUDDdistribute2}, when applied to all the balanced runs of a rotor, gives the inequality

\begin{equation}\label{balancecoeff}
b(UU(r)) + b(DD(r)) \geq 2b(r).
\end{equation}

This can be seen from the fact that the balanced runs in two periods of $r$ are split between $UU(r)$ and $DD(r)$. Equality occurs iff no balanced run is split up as in Theorem \ref{UUDDdistribute2}. Next we define the analogous decomposition of a rotor into uniform runs.

\begin{defn}
Define an \textit{uniform run decomposition} (URD) of a rotor $r$ as a partition of the period of $r$ into contiguous uniform runs. Define the uniform run type of this URD as the infinite periodic sequence $(u_1, u_2, \ldots)$ such that $u_i$ is the length of the $i$-th run in this URD, $i$ taken modulo the number of runs in the URD.
\end{defn}

A rotor can have multiple URD's if we choose to split up larger uniform runs. Thus $122112$ can have URD $1|22|11|2$ or $1|2|2|11|2$ or $1|2|2|1|1|2$. Based on these two definitions and inequality (\ref{balancecoeff}), we can conclude the following:

\begin{thm}\label{BURD}
If $r \in \mathbf{BAL}$ has a BRD of type $(b_1, b_2, \ldots)$, then one of the following must hold:

\begin{enumerate}[i.]
\item One of $b(UU(r))$ and $b(DD(r))$ is strictly larger than $b(r)$.
\item There is also a URD of $r$ of type $(b_1, b_2, \ldots)$.
\end{enumerate}

\begin{proof}
The theorem follows from the fact that either a split of the type given in Theorem \ref{UUDDdistribute2} exists for $UU(r)$ or $DD(r)$, in which case $b(UU(r)) + b(DD(r)) > 2b(r)$, since an extra balanced run was created, or else every balanced run of length $2b_i$ corresponds to a uniform run of length $b_i$ at the position half its own.
\end{proof}
\end{thm}

We define a \textit{BURD rotor} as a rotor type satisfying condition \textit{ii} of Theorem \ref{BURD}. That is, it has a BRD and a URD of the same type. From any balanced rotor, repeated applications of the compressor must eventually lead to a BURD rotor; otherwise we could increase the balance coefficient arbitrarily often using $UU$ or $DD$, and it can only take on a finite number of rational values. Naturally, the class of BURD rotors is denoted $\mathbf{BURD}$. One essential property of BURD rotors is the following generalization of Corollary \ref{abbafixed}:

\begin{thm}\label{BURDfixed}
If $r\in \mathbf{BURD}$, then $UD(r) \equiv DU(r) \equiv r$.

\begin{proof}
Assume we are looking at $UD(r)$, without loss of generality. Let a uniform run of the URD of $r$ that has a BRD counterpart be $r^{(a + 1)}, r^{(a + 2)}, \ldots, r^{(a+b)}$. The corresponding balanced run in the BRD is $r^{(2a + 1)}, r^{(2a+2)}, \ldots, r^{(2a+2b)}$. Thus, when the particles numbered $(2a+1)$ through $(2a+2b)$ are emitted by the source rotor, exactly half go to rotor 2 of the compressor, hitting the uniform run corresponding to $r^{(a+1)}, \ldots, r^{(a+b)}$ there, and the other half go to the same uniform run in rotor 3. Only one of these runs is routed down, depending on $r^{(a+1)}$, so we know that this uniform run is copied faithfully onto $UD(r)$. Each of the uniform runs in the URD is copied in the same way into $UD(r)$. It is now evident that $UD(r) \equiv r$, and the proof of $DU(r) \equiv r$ is exactly analogous.
\end{proof}
\end{thm}

Note that the fixed cycles of the compressor are the only rotors necessary to consider to prove Conjecture \ref{central}. A \textit{fixed cycle} of the compressor is defined formally as a set $\mathbf{C}$ of rotors, closed under all compressor operations, such that any element can be reached from any other by way of some sequence of compressor operations. From the preceding theorems, a standard extremal argument shows:

\begin{thm}\label{fixedcycle}
For every fixed cycle $\mathbf{C}$ of the compressor $\mathbf{C} \subset \mathbf{BURD}$.
\begin{proof}
Suppose for the sake of contradiction that some such cycle $\mathbf{C}$ exists, containing a rotor $r \not\in \mathbf{BURD}$. We can choose $r$ to have maximal $b(r)$ of all rotors in $\mathbf{C} \setminus \mathbf{BURD}$. Also, $\mathbf{C} \cap \mathbf{BURD}$ is nonempty, so we can choose $r' \in \mathbf{C} \cap \mathbf{BURD}$ with maximal balance coefficient in this set as well. Clearly, $b(r') > b(r)$ since there is a sequence of $UU$'s and $DD$'s strictly increasing $b(r)$ until it becomes BURD, by Theorem \ref{BURD}. But then, if $r'$ can reach $r$, it must also be able to reach some rotor $r''$ with $b(r'') > b(r')$ by inequality \ref{balancecoeff}, a contradiction. Thus $\mathbf{C} \subset \mathbf{BURD}$ always.
\end{proof}
\end{thm}

The obvious class of BURD rotors is $\mathbf{ABBA}$, and the generalization thereof to block-repetitive rotors (i.e. aabb-bbaa rotors and the like). However, not all BURD rotors are of this form. For instance, the rotor type $1122222111$ generates a nontrivial fixed cycle of size $2$ in the action of the compressor. However, when the BURD rotor also happens to be ab-ba, we have some further results.

\section{On ab-ba rotors}\label{abba}

Recall that a balanced rotor $r$ is called ab-ba if it is $2$-balanced, that is, if $r^{2k+1}$ and $r^{2k+2}$ are different for all $k \in \mathbb{N}$. Intuitively, ab-ba rotors are just one step away from the universal rotor $12$. This is because the compressor keeps ab-ba rotors ab-ba, and they are all fixed points of $UD$ and $DU$. Also, they can't be block-repetitive. Hence, to be boppy, they need to be palindromic and thus have length $0 \pmod{4}$ (consider the middle two terms). If our conjectures hold, then all ab-ba rotors of length $2 \pmod{4}$ should be universal. Indeed, we verified using computer programs that all ab-ba rotors of length $2\pmod{4}$ up to length $42$ are universal by the compressor algorithm. Now we discuss some extra tools.

\begin{defn}
Break an ab-ba rotor $r$ into its constituent $12$ or $21$ blocks. If exactly $k$ such blocks are $21$, then define the \textit{ba-frequency} of $r$ as $\displaystyle m(r) = \frac{2k}{|r|}$.
\end{defn}

This function $m$ can be made into a monovariant for ab-ba rotors under the compressor, as in the following result.

\begin{thm}\label{bafreq}
For $r\in \mathbf{ABBA}$, the inequality $m(UU(r)) + m(DD(r)) = 2m(r)$ holds.

\begin{proof}
We only need the number of $21$'s in the first $|r|$ terms of $UU(r)$ and the number of $21$'s in the first $|r|$ terms of $DD(r)$ to sum to the number of $21$'s in the first $2|r|$ terms of $r$. Each $21$ in the first two periods of $r$ corresponds to one $21$ in the first $|r|$ terms of either $UU(r)$ or $DD(r)$. For instance, suppose that $r^{(k)} = 1$, $r^{(2k-1)} = 2$, and $r^{(2k)} = 1$. Then, in configuration $UU$, $($rotor $2)^{(k)}$ and $($rotor $3)^{(k)}$ both go back to the source, while in $DD(r)$, they route to targets $4$ and $5$. In $UU(r)$ this block $21$ will not appear in the hitting sequence, while in $DD(r)$ it contributes a $21$. Each block $21$ in the $2|r|$ terms of $r$ contributes a $21$ to exactly one of $UU(r)$ and $DD(r)$.
\end{proof}
\end{thm}

In particular, this shows that there is some sequence of $UU$'s and $DD$'s that never increase $m(r)$. Further, the theorem shows that if $m(r)$ is changed at all by one of $UU$ and $DD$, then it must be strictly decreased by one of them. We end with two specific results based on this idea.
\begin{thm}\label{one21}
Every $r\in \mathbf{ABBA}$ consisting of more than one block of $12$ and exactly one block of $21$ is universal.

\begin{proof}
Let $f(r)$ denote the position of the $21$ block in such an $r \in \mathbf{ABBA}$. If either of $m(UU(r))$, $m(DD(r))$ is not equal to $m(r)$, then one of them must be strictly less than $m(r)$ and be zero, by Theorem \ref{bafreq}. Otherwise, both $UU(r)$ and $DD(r)$ also have exactly one $21$ block.

If $|r| \equiv 0 \pmod{4}$, then one of $UU(r)$ and $DD(r)$ is $12$ unless $f(r) = 1$ or $f(r) = \frac{|r|}{2}$. If $f(r) = 1$, then one of $f(UU(r))$, $f(DD(r))$ is $\frac{|r|}{4} + 1$, and another application of $UU$ or $DD$ suffices. If $f(r) = \frac{|r|}{2}$, then one of $f(UU(r))$, $f(DD(r))$ is $\frac{|r|}{4}$, and again $r$ is universal after another application of $UU$ or $DD$.

If $|r| \equiv 2 \pmod{4}$, then one of $UU(r)$ and $DD(r)$ is $12$ if $f(r) = 1$ or $f(r) = \frac{|r|}{2}$. Otherwise, it is easy to prove that one of $f(UU(r))$, $f(DD(r))$ is $\displaystyle \left\lceil \frac{f(r)}{2}\right\rceil$, so some combination of $UU(r)$, $DD(r)$ must move the $21$ block to position $1$.
\end{proof}
\end{thm}

The proof above shows that the ab-ba problem can be reinterpreted as follows: \vspace{5mm}

Given a periodic sequence $A = \{a_i\}_{i = 0}^{\infty}$ of $0$'s and $1$'s, where $0$ represents $12$ and $1$ represents $21$, we are allowed to perform the following two transforms on the sequence. The first takes $A$ to $UU(A) = \{b_i\}_{i = 0}^{\infty}$, and the second to $DD(A) = \{c_i\}_{i = 0}^{\infty}$, where $b_i = a_{2i - a_i}$ and $c_i = a_{2i - 1 + a_i}$. The problem is to determine when these two transforms, applied repeatedly, can take a sequence to a constant one. In this spirit, the next theorem is immediate.

\begin{thm}
If $r\in \mathbf{ABBA}$ satisfies the following conditions:
\begin{enumerate}[i.]
\setlength\itemindent{25pt} 
\item $|r| \equiv 0 \pmod{4}$.
\item For any two $21$ blocks of $r$ that are the $a$-th and $b$-th blocks of $r$, $a$ and $b$ do not satisfy $\lceil \frac{a}{2}\rceil \equiv b \pmod{|r|}$.
\item For every pair of $21$ blocks of $r$, their positions as blocks have the same parity.
\end{enumerate}

then $r\in \mathbf{U}$.

\begin{proof}
It is clear from the proof of Theorem \ref{one21} that one of $UU(r), DD(r)$ is $12$, since each $21$ block acts independently of the others.
\end{proof}
\end{thm}

If we are patient we can extend these two results somewhat, but it does not seem fruitful in general to do so.

\section{Directions for Future Research}\label{future}
Here we collect together some possible directions of future work on the universality problem and the related string manipulations. The algebraic perspective in Section \ref{compressor} may have further implications for the compressor algorithm. In Section \ref{empiricalresults} we mentioned several conjectures inspired by empirical data. It was conjectured that the compressor directly proves the universality of all non-palindromic rotors of prime length and all ab-ba rotors of length $2\pmod{4}$; these are the focus of our current work. Another direction to take is the investigation of the universality of BURD rotors, which are central to the compressor algorithm. The classification of BURD strings is interesting in its own right, especially given that nontrivial examples like 1122222111 exist. Also, the fascinating self-reference of the BURD condition merits more attention.

Another problem proposed by Propp asks for the universality of rotors when \textit{cyclic shifts of the same rotor are considered equivalent}. Our results point to the conjecture that all rotors are universal under this alternative definition. The universality question can also be generalized in another direction, that of determining which finite sets of rotors are universal when taken together. Resolution of the universality problem would almost completely solve this as well. We may try to quantify our results on universality; instead of asking whether or not rotor $r$ makes rotor $r'$, we can ask \textit{how efficiently} $r$ can make $r'$. Such a question provides another notion of the ``complexity" or ``information value" of binary strings.

\section{Acknowledgements}\label{acknowledgements}

The author is grateful to Tanya Khovanova of MIT for her invaluable guidance and insight throughout this project. She provided much-needed advice during the writing of this paper. The author also received significant support from James Propp of U. Mass. Lowell, who proposed the problem and provided expertise on rotor-routers. The present work was supported by the MIT-PRIMES program for high school students in the Greater Boston area. Thanks go out to Yan Zhang of MIT for providing a philosophical comment on the responsibility of memory (to be found in Section \ref{examples}), and David Yang of Philips Exeter Academy for proofreading this paper.

\bibliographystyle{alpha}

\end{document}